\documentclass{cicp}  

\usepackage{amsmath,amsfonts}
\usepackage{algorithm}
\usepackage{array}
\usepackage[caption=false,font=normalsize,labelfont=sf,textfont=sf]{subfig}
\usepackage{textcomp}
\usepackage{stfloats}
\usepackage{url}
\usepackage[table]{xcolor}
\usepackage{tabularx}
\usepackage{booktabs}
\usepackage{pifont}
\usepackage{verbatim}
\usepackage{latexsym}
\usepackage{graphicx}
\usepackage{cite}
\usepackage[hidelinks]{hyperref}
\usepackage{algpseudocode}
\usepackage{subcaption}
\usepackage{stmaryrd}
\usepackage{tikz}

\begin{document}
	
	\title{Learning Stochastic Hamiltonian Systems via Stochastic Generating Function Neural Network}
	
	\author[Chen C et.~al.]{Chen Chen\affil{1}, Lijin Wang\affil{1}\corrauth, Yanzhao Cao\affil{2}, Xupeng Cheng\affil{1}}
	\address{\affil{1} \ School of Mathematical Sciences, University of Chinese Academy of Sciences, Beijing 100049, P.R. China.\\
		\affil{2} \ Department of Mathematics and Statistics, Auburn University, Auburn, AL 36849, USA.}
	\email{{\tt ljwang@ucas.ac.cn} (L.J.Wang)}
	
	
	\begin{abstract}
		In this paper we propose a novel neural network model for learning stochastic Hamiltonian systems (SHSs) from observational data, termed the stochastic generating function neural network (SGFNN). SGFNN preserves symplectic structure of the underlying stochastic Hamiltonian system and produces  symplectic predictions. Our model utilizes the autoencoder framework to identify the randomness of the latent system by the encoder network, and detects the stochastic generating function of the system through the decoder network based on the random variables extracted from the encoder. Symplectic predictions can then be generated by the stochastic generating function. Numerical experiments are performed on several stochastic Hamiltonian systems, varying from additive to multiplicative, and from separable to non-separable SHSs with single or multiple noises. Compared with the benchmark stochastic flow map learning (sFML) neural network, our SGFNN model exhibits higher accuracy across various prediction metrics, especially in long-term predictions, with the property of maintaining the symplectic structure of the underlying SHSs.
		
	\end{abstract}
	
	
	
	\maketitle
	\section{Introduction}
	\label{sec1}
	Hamiltonian systems describe the motion of dynamical systems where dissipation is negligible, such as in the fields of celestial mechanics, aerospace engineering and life science etc.  The essential characteristic of Hamiltonian systems is that their phase flows preserve the symplectic structure, which can be geometrically interpreted as area preservation in the phase space (\cite{Feng,Hairer}). With the discovery of more and more random features in real world and the development of uncertainty quantification, growing attentions are drawn to stochastic modeling in science and technology, leading to increasing interest in the study of stochastic Hamiltonian systems, which are Hamiltonian systems perturbed by noises while still maintaining the symplectic structure almost surely in the phase space.
	Bismut et al.\cite{Bismut} investigated the fundamental theory of stochastic Hamiltonian systems and explored how stochastic perturbations affect the system dynamics. Misawa \cite{Misawa} discussed conservation laws and symmetry issues in stochastic dynamical systems. Milstein, Repin, and Tretyakov \cite{Milstein_mul,Milstein_add} proposed symplectic integration methods for Hamiltonian systems with additive and multiplicative noises, initiating systematic research on symplectic algorithms for stochastic Hamiltonian systems. Subsequently, there arise a series of research on stochastic symplectic numerical methods, including \cite{Hong,Wang,Bou, Deng} etc, for which more references and detailed introduction can be found in the monograph \cite{HongSun}.
	
	In recent years, with the development of data science and machine learning, detecting governing differential equations of dynamical systems from observational data via neural networks has been absorbing more and more attention. There has been much work devoting to the learning of ordinary differential equations (ODEs)\cite{Node, Zhong}, partial differential equations (PDEs)\cite{Long1, Long2} and stochastic differential equations (SDEs)\cite{Kong, Nsde}. In \cite{HNN}, Greydanus et al. proposed to set up neural networks that respect the physical conservational laws of the latent dynamical systems, and constructed the Hamiltonian neural network (HNN) which learns the Hamiltonian function of a Hamiltonian system instead of its total vector field, unlike usual ODE-nets, and can preserve the energy (Hamiltonian function) of the system. This opens up the structure-preserving neural network learning of Hamiltonian systems. 
	
	From then on, a number of studies have been devoted to the symplectic learning of Hamiltonian systems, aiming to preserve the intrinsic symplecticity during neural network-based learning and prediction. A notable example is the Symplectic Recurrent Neural Network (SRNN) introduced by Chen et al.~\cite{SRNN}, which extends the idea of the HNN by incorporating the leapfrog integrator as a dynamic generator. By modeling the Hamiltonian as a sum of two networks, SRNNs effectively learn and predict separable Hamiltonian systems while maintaining symplectic structure throughout both training and inference. Zhu et al.~\cite{Zhu} advanced the HNN's line of work by employing the inverse modified equation of a symplectic method to compute the Hamiltonian function. Their results established that HNNs built upon symplectic schemes possess well-defined learning targets, thereby providing theoretical validation for their superior performance. Complementing this, Jin et al.~\cite{Jin} proposed SympNets, a modular architecture composed of simple symplectic blocks, and showed that given suitable activation functions, these networks can approximate arbitrary symplectic maps directly from data. To incorporate additional structural priors, Tong et al.~\cite{Tong} embedded Taylor series expansions into neural networks to enforce symmetry, specifically targeting the learning of Hamiltonian gradients for separable systems. Their approach utilized a fourth-order symplectic Runge-Kutta method to generate the associated symplectic maps. Toth et al.~\cite{Toth} took a generative modeling perspective and developed Hamiltonian Generative Networks, which extract low-dimensional abstract states from high-dimensional observational data, learn separable Hamiltonian dynamics in this latent space, and reconstruct the high-dimensional trajectories using a decoder network. In a related approach, the work in~\cite{GFNN} introduced a fully connected deep neural network trained to approximate the generating functions of deterministic Hamiltonian systems. This method directly learns the symplectic map from data and guarantees that the global prediction error grows at most linearly over time. Other contributions to this growing field include the models proposed by Xiong et al.~\cite{Xiong}, Wu et al.~\cite{Wu}, and Dierkes et al.~\cite{Dierkes}, each offering new insights and architectures for learning deterministic Hamiltonian systems through neural networks.

	As far as data-based neural network inversion of stochastic differential equations (SDEs) is concerned, the essential challenge is to quantify and model the randomness underlying the data. 
	Kong et al. \cite{Kong} proposed the SDE-Net, a method that treats the transformations of deep neural networks as evolutions of stochastic dynamical systems, which uses Drift Net and Diffusion Net to learn the deterministic and uncertain parts of the system, respectively, providing a new approach to quantify uncertainty in deep neural networks. Dietrich et al. \cite{Dietrich} utilized the neural network to learn the drift and diffusion coefficients of SDEs. They designed the loss function by incorporating classical stochastic numerical integrators. 
	The challenge introduced by randomness fundamentally stems from the unobservability of the noise in the data. To address this issue, various denoising techniques have been adopted in the literature. Many methods leverage the evolution of the probability density of the system's state for inference, such as learning SDEs based on the Fokker-Planck equation \cite{FK1, FK2}, variational methods \cite{variational1, variational2}, maximum likelihood estimation \cite{Dietrich}, normalizing flows \cite{normalizing1, normalizing2}. Chen et al. \cite{Chen1, Chen2, Chen3}, Xu et al. \cite{Xu} extended deterministic flow map methods to stochastic systems, using generative stochastic models to learn stochastic systems. Liu et al. \cite{Liu1, Liu2} developed a training-free conditional diffusion model to learn unknown SDEs from observational data. This method generates labels by solving inverse-time ODEs or inverse-time SDEs, enabling supervised learning of the flow map. 
	
	For the structure-preserving learning of stochastic Hamiltonian systems, there has been limited research effort, with notable contributions from \cite{WangZhan1}, \cite{WangZhan2}, and \cite{Cheng}. These studies primarily extended the idea of HNN to stochastic Hamiltonian systems, involving learning the drift and diffusion Hamiltonian functions to capture the symplectic structure of the system, thereby enabling symplectic predictions. However, the Hamiltonian functions 
	are only approximate generators of the symplectic solution flow, i.e. there is discrepancy between the symplectic mappings generated by the Hamiltonian functions via symplectic schemes and the true solution flow mappings, for which the works \cite{Zhu,france} have resorted different modified equations as remedy. This motivates our effort to learn directly the exact stochastic generating function of the symplectic solution flow mappings for the SHSs. 
	
	In this paper, we propose a neural network method for extracting the stochastic generating function of unknown stochastic Hamiltonian systems from observational data. Our model is based on the autoencoder architecture. The encoder component is utilized to identify Gaussian hidden variables in adjacent data pairs, and the decoder part employs the learned hidden variables to reconstruct the stochastic generating function. Given that the stochastic generating function serves as a random function intrinsically linked to stochastic Hamiltonian systems and can implicitly define their exact symplectic flow maps\cite{Wang,Deng}, our approach enables the accurate representation of the underlying symplectic map embedded within the observational data. To ensure the effectiveness of our model, we carefully design the loss function to constrain the latent variables to follow a standard Gaussian distribution and to guarantee that the decoder's output closely approximates the true stochastic generating function. Consequently, the reconstructed stochastic generating function not only captures the dynamical characteristics of the system but also preserves its geometric structure, leading to more reliable long-term predictions and simulations.
	We validate our approach through numerical experiments on four distinct stochastic Hamiltonian systems, using the sFML model \cite{Xu} as a baseline for comparison. The experimental results demonstrate that our model achieves significantly improved 
	learning effect than the baseline model, in terms of the accuracy of the predicted probability density, mean and standard deviation of the solution, as well as in reflecting the geometric or physical feature of the underlying dynamical system, such as the conservation of the Hamiltonians and the linear growth property of the solution's second moments, especially in long-term situations. 

	The rest of the paper is organized as follows. Section \ref{sec2} introduces the stochastic Hamiltonian systems and the learning objectives. Section \ref{sec3} presents the detailed construction of the stochastic generating function neural network, including the network architecture and the loss function. In section \ref{sec4}, we conduct numerical experiments on learning various stochastic Hamiltonian systems with our model, including separable and non-separable, additive and multiplicative stochastic Hamiltonian systems with single or multiple noises, compared with the benchmark sFML model. Numerical results demonstrate the effectiveness and superiority of our model. Section \ref{sec5} is a brief conclusion.
	
	\section{Problem statement}
	\label{sec2}
	A stochastic Hamiltonian system refers to a Hamiltonian system subject to certain random perturbations. We consider an autonomous 2d-dimensional stochastic Hamiltonian system of Stratonovich type(\cite{Milstein_mul,Milstein_add}).
	
	\begin{equation}\label{eqn1} 
		\begin{aligned}
			d p=f(p, q) d t+\sum_{k=1}^r \sigma_k(p, q) \circ d W_k(t), & \quad p\left(t_0\right)=p_0, \\
			d q=g(p, q) d t+\sum_{k=1}^r \gamma_k(p, q) \circ d W_k(t), & \quad q\left(t_0\right)=q_0,
		\end{aligned}
	\end{equation}
	for which there exist smooth functions $H(p,q),\,\,H_k(p,q), \,\,k=1,...,r$ such that
	\begin{equation}\label{eqn2}
		\begin{aligned}
			& f(p, q)=-\frac{\partial H}{\partial q}(p, q)^T,\quad\sigma_k(p, q)=-\frac{\partial H_k}{\partial q}(p, q)^T,\\
			& g(p, q)=\frac{\partial H}{\partial p}(p, q)^T,\quad\gamma_k(p, q)=\frac{\partial H_k}{\partial p}(p, q)^T,
		\end{aligned}
	\end{equation}
	where \(p, q, f, g, \sigma_k, \gamma_k\) are \(d\)-dimensional vectors. \(W_k(t)\), \(k = 1, 2, \ldots, r\), are independent standard Brownian motions defined on a complete filtered probability space \(\{ \Omega, \{ \mathcal{F}_t \}_{t \geq 0}, \mathbb{P} \}\). 
	
	We assume that \(f\), \(g\), \(\sigma_k\) and \(\gamma_k\) are unknown, and the Brownian motions \(W_k\) are unobservable, k=1,...,r. We only observe the discrete-time data of the solution flow \((p, q)\) at time points within the interval [0,T]:
	\begin{equation}
		\{(p_{t_0}^{(i)}, q_{t_0}^{(i)}), (p_{t_1}^{(i)}, q_{t_1}^{(i)}), \ldots, (p_{t_{L_i}}^{(i)}, q_{t_{L_i}}^{(i)})\}, \quad i = 1, 2, \ldots, N.
	\end{equation}
	For simplicity, we make the following assumptions: the length of the observed trajectories are equal, i.e., \(l_i = L\), i=1,...,N, and the intervals between discrete time points are also equal, i.e., \(t_{j+1} - t_j = \Delta\), j=0,...L-1. Therefore, our observational data have the following form:
	\begin{equation}\label{data}
		\{x_{0}^{(i)}, x_{1}^{(i)}, \ldots, x_{L}^{(i)}\}, \quad x_{j}^{(i)} = (p_{t_j}^{(i)T}, q_{t_j}^{(i)T})^T, \quad i = 1, 2, \ldots, N, \quad j = 0, 1, \ldots, L.
	\end{equation}
	
	Let $x_t:=(p_t^T,q_t^T)^T$. The solution flow $\varphi_t(x_0; \omega): x_0\to x_{t,\omega}$ of a SHS preserves the symplectic structure almost surely (\cite{Milstein_add}), where \(\omega \in \Omega\) denotes a sample point in the sample space $\Omega$, that is 
	
	\begin{equation}\label{eqn3}
		\begin{aligned}
			&{({\varphi_{t_0}}')}^T J {\varphi_{t_0}}'=J, \quad \mbox{almost surely}, \\
		\end{aligned}
	\end{equation}
	where $\varphi_{t_0}':=\frac{\partial \varphi_t}{\partial x_0}$, $J=\begin{bmatrix} 0 & I \\ -I & 0 \end{bmatrix}$ is an $2d \times 2d$ matrix, and $I$ is a $d \times d$ identity matrix.
	
	Our goal is to learn the solution flow \(\varphi_t\) of SHS (\ref{eqn1}) using observational data (\ref{data}). Assume the learned map is $ \widetilde{\varphi} \approx \varphi_{\Delta}$. Given the initial point $ x_0$, we can use $ \widetilde{\varphi} $ to generate a discrete trajectory that approximates the true solution. We achieve this by using a mathematical tool called stochastic generating function \cite{Wang,Deng}, which ensures that the symplectic property of \(\widetilde{\varphi}\) is preserved. Our method involves using an autoencoder network model to approximate the stochastic generating function. Once the stochastic generating function is determined, the symplectic map is also learned.
	
	\section{Stochastic generating function neural network}
	\label{sec3}
	
	
	In this section, we propose a Stochastic Generating Function Neural Network (SGFNN) based on the autoencoder framework to learn the stochastic generating function of SHS. 
	
	A standard autoencoder network consists of an encoder that compresses the input into a latent space and a decoder that reconstructs the original input from this latent representation. In \cite{Xu}, the autoencoder network is used to implement the stochastic flow map learning (sFML) method, which directly learns the mapping \(x_{n+1} = G_\Delta(x_n, z_n)\) on the SDE solution trajectory from the state $x_n$ at the time $t_n$ to the state $x_{n+1}$ at the time $t_{n+1}=t_n+\Delta$, where $z_n$ drawn from a unit Gaussian distribution governs the stochasticity of the mapping. 
	
	However, for the learning of stochastic Hamilton systems, the sFML method can not naturally preserve the symplecticity of the flow map. Since the stochastic generating functions can generate the symplectic flows of SHSs \cite{Wang,Deng}, we propose to learn the stochastic generating function via the autoencoder network, and then use it to generate the symplectic flow map of the latent SHS. 
	
	In \cite{GFNN}, GFNN method is proposed to learn the generating functions of deterministic Hamiltonian systems. Our SGFNN is designed for stochastic Hamiltonian systems which employs the autoencoder network to deal with the stochasticity.
	The training data of the SGFNN are pairs of states $\{(x_0^{(i)},x_{1}^{(i)})|i=1,\dots,M\}$ from evolution trajectories of a SDE with time interval $\Delta$. The encoder  extracts the latent random variable $ z^{(i)}$ underlying each data pair $(x_0^{(i)},x_1^{(i)})$ for $i=1,\dots,M$, and the decoder reconstructs the stochastic generating function $ S $. 

	\subsection{Stochastic generating function}

		Stochastic generating functions can generate symplectic mappings locally, and generate the symplectic flow of stochastic Hamiltonian systems almost surely \cite{Deng,Wang}. There are different types of generating functions depending on different coordinates \cite{Hairer,Deng,Wang,Hong}, and we use here the first kind which generates symplectic mappings according to the following theorem.
		\begin{theorem}[{\cite{Deng}}]
			Let $\varphi_{\omega} :(p,q)\rightarrow (P,Q)$ be a smooth random map from $\mathbb{R}^{2d}$ to $\mathbb{R}^{2d}$ with $\omega$ being a sample point from the underlying probability space. Then $\varphi_{\omega}$ is symplectic if there exists locally a smooth random map $S(P,q,\omega)$ from $\mathbb{R}^{2d}$ to $\mathbb{R}^{2d}$ such that $\frac{\partial(P^T q+S(P,q,\omega))}{\partial P \partial q}$ is invertible a.s. and we have a.s.
			\begin{align}\label{generating1}
				p=P+\frac{\partial S}{\partial q}^T(P,q,\omega),\quad Q=q+\frac{\partial S}{\partial P}^T(P,q,\omega).
			\end{align}
		\end{theorem}
		The function $S$ is called the stochastic generating function (of the first kind), which generates a symplectic mapping $(p,q)\rightarrow (P,Q)$ according to the relation \eqref{generating1} under the almost sure invertibility of $\frac{\partial(P^T q+S(P,q,\omega))}{\partial P \partial q}$.
		

		\subsection{Method description}
		Given the observational data pairs $\{(x_0^{(i)},x_1^{(i)})|i=1,\dots,M\}$ of a latent stochastic Hamiltonian system, which can be regarded as being generated by a stochastic generating function $S(p_1,q_0,z)$ according to the relation \eqref{generating1}. Then we attempt to utilize an autoencoder network model, where the encoder captures the latent random variable \(z\) which follows a unit Gaussian distribution. The decoder outputs the generating function \(S_{\theta}(p_1, q_0, z)\) with $\theta$ representing the learnable parameters in the neural network. The final loss function is composed of two parts: one part constrains the output of the encoder \(z\) to follow a unit Gaussian distribution, and the other part ensures that the output of the decoder \(S_{\theta}\) approximates the true \(S\).
		
		\subsubsection{Encoder part}
		The time interval between $ x_0^{(i)} $ and $ x_1^{(i)} $ is $ \Delta $.
		The encoder outputs the latent variable \(z^{(i)}\) for each input observational data pair \((x_0^{(i)}, x_1^{(i)})\). 
		
		The latent random variable \(z\) needs to satisfy two constraints: first, it must be independent of \(x_0\); second, it should follow a unit Gaussian distribution. Note that \(z\) can be a multidimensional random variable. Since these two conditions are statistical constraints on \(z\), we input $ K $ data pairs at once each time during the training, for which the encoder outputs $ K $ latent variables $ z $ to form the set $B$, i.e.
		\begin{equation}
			B = \{ z^{(i)} = E_\Delta(x_0^{(i)}, x_1^{(i)},\theta_1)|i=1,\dots,K \},
		\end{equation}
		where $\theta_1$ are network parameters of $E_{\Delta}$.
		
		The training set contains \(N_B\) batches, and each batch consists of \(K\) data pairs. We use the sub-sampling strategy described in \cite{Xu} to guarantee that the distribution of \(z\) is independent of \(x_0\). This strategy is achieved by adjusting the data in the training set to ensure that the distribution of $ \{{x_0^{(i)}}\}_{i=1}^{K} $ is non-Gaussian for all \(N_B\) batches. For the $ K $ data pairs $\{(x_0^{(i)}, x_1^{(i)})| i=1,\ldots,K\}$ input to the encoder, $ \{{x_0^{(i)}}\}_{i=1}^{K} $ represents the set of $ K $ points closest to a given $ x_0 $ that is randomly chosen from the observed data. 
		
		In the encoder part, the distribution loss function is denoted as \(L_D(B)\), which is used to constrain the probability distribution of \(z\) to approximate a unit Gaussian distribution. This distribution loss function is constructed in the same way as \cite{Xu}:
		\begin{equation}
			\mathcal{L}_D(B) = \mathcal{L}_{\text{distance}}(B) + \tau \cdot \mathcal{L}_{\text{moment}}(B),
		\end{equation}
		where \(\mathcal{L}_{\text{distance}}(B)\) measures the distance between the distribution of \(z\) in \(B\) and the unit Gaussian distribution, and \(\tau\) is a hyperparameter that balances the two terms. Concretely, 
		\begin{equation}
			\mathcal{L}_{\text{distance}}(B) = \|\hat{f}_B - f_{normal}\|_2,
		\end{equation}
		where \(\hat{f}_B\) is the probability density function of $ z $ obtained using the kernel density estimation (KDE) method, and \(f_{normal}\) is the density function of unit Gaussian distribution. \(\mathcal{L}_{\text{moment}}(B)\) measures the discrepancy between the central moments of \(z\) in \(B\) and those of the standard Gaussian distribution. Considering that the dimension $n_z$ of the latent variable \(z\) may be greater than 1, \cite{Xu} adds the cross-correlation coefficients  $\widehat{\rho}_{j k}(B)$ between the elements of \(z\) in \(\mathcal{L}_{\text{moment}}\).
		\begin{equation}
			\mathcal{L}_{\text {moment }}(B)=\sum_{j=1}^6 \frac{1}{c_j}\left\|\widehat{\mu}^{(j)}(B)-\mu^{(j)}\left(P_{normal}\right)\right\|_2^2+\frac{\nu}{V} \sum_{1 \leq j<k \leq n_z}\left(\widehat{\rho}_{j k}(B)\right)^2,
		\end{equation}
		where \(\hat{\mu}^{(j)}(B)\) denotes the \(j\)-th central moment of \(z\) in the batch \(B\), and \(\mu^{(j)}(P_{normal})\) denotes the \(j\)-th central moment of unit Gaussian. $
		{V}=n_z\left(n_z-1\right) / 2, \quad v=0.1, \quad 
		c_1=1.0, \quad c_2=1.0, \quad c_3=2.0, \quad c_4=3.0, \quad c_5=8.0, \quad c_6=15.0 .$ 
		
		\subsubsection{Decoder part}
		In the decoder part, we aim to approximate the stochastic generating function $S$. The latent variable \(z\), which follows a unit Gaussian distribution, has already been identified by the encoder. According to equation \eqref{generating1}, \(S\) is a function of \((p_1, q_0,z)\), therefore the input to the decoder includes three sets of data: the value of \(p_1\) at the time $t_0+\Delta$, the value of \(q_0\) at the time $t_0$, and the latent variable \(z\), and 
		\begin{equation}
			S_{\theta_2} = D_\Delta(p_1^{(i)}, q_0^{(i)}, z^{(i)},\theta_2), \quad i = 1, \ldots, M,
		\end{equation}
		with learnable network parameters $\theta_2$.
		According to equation (\ref{generating1}), we define the mean squared error (MSE) loss on $B$ as:
		\begin{equation}
			\mathcal{L}_{M S E}(B)=\frac{1}{K} \sum_{i=1}^K\left(\left\|p_1^{(i)}-p_0^{(i)}+\frac{\partial S_{\theta_2}}{\partial q_0^{(i)}}^T\right\|_2^2+\left\|q_1^{(i)}-q_0^{(i)}-\frac{\partial S_{\theta_2}}{\partial p_1^{(i)}}^T\right\|_2^2\right),
		\end{equation}
		where $\frac{\partial S_{\theta_2}}{\partial q_0^{(i)}}$ and $\frac{\partial S_{\theta_2}}{\partial p_1^{(i)}}$ denotes the partial derivative of $S_{\theta_2}$ with respect to $q_0$ and $p_1$ respectively, which are evaluated at the point $(p_1^{(i)},q_0^{(i)},z^{(i)})$. The total training loss is:
		\begin{equation}
			\mathcal{L}_{SGFNN}=\frac{1}{N_B} \sum_{j=1}^{N_B} \mathcal{L}\left(B_j\right)=\frac{1}{N_B} \sum_{j=1}^{N_B}\left[\mathcal{L}_{MSE}\left(B_j\right)+\lambda \cdot \mathcal{L}_D\left(B_j\right)\right],
		\end{equation}
		where \(\lambda\) is a hyperparameter that balances  $\mathcal{L}_{MSE}$ and $\mathcal{L}_D$.
		
		
		\begin{figure}[!t]
			\centering
			\includegraphics[width=0.60\textwidth, bb=0 0 258 113]{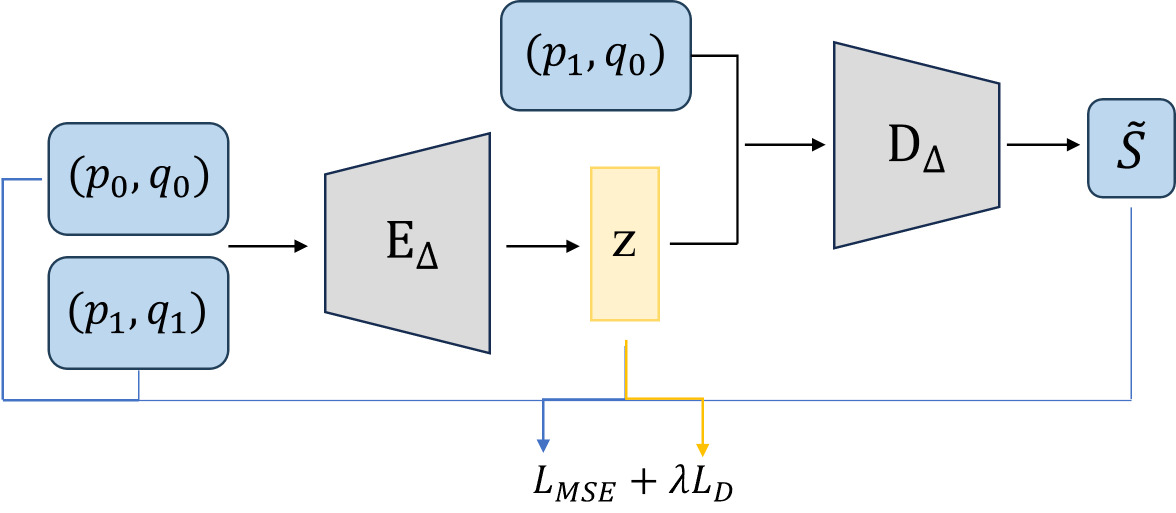} 
			\caption{An illustration of the algorithm flow and loss function of the SGFNN Network} 
			\label{figtrainingprocess} 
		\end{figure}

		The complete training process of the SGFNN is illustrated in Figure~\ref{figtrainingprocess}. Minimizing \(\mathcal{L}_D\) aims to adjust the parameters of the encoder network so that the output \(z\) matches the desired distribution. Minimizing \(\mathcal{L}_{\text{MSE}}\) aims to adjust the parameters of both the encoder and decoder networks so that the output \(S_{\theta_2^*}\) approximates the true stochastic generating function \(S\).

		\subsubsection{Testing stage}
		After training, we use only the decoder part for prediction. Given the system data \(x_0=(p_0^T, q_0^T)^T\) at the current instant $t=t_0$ and a random variable \(\omega\) that follows a unit Gaussian distribution, we need to predict the \(\widetilde{x}_1=(\widetilde{p}_1^T, \widetilde{q}_1^T)^T\) at the next instant $t=t_0+\Delta$ . Since the symplectic map \((p_0, q_0) \to (p_1, q_1)\) defined by equation (\ref{generating1}) is implicit, we use the fixed-point iteration to compute the data of the next instant, which is
		\begin{align}
			p_1^{j+1} &= p_0 - \frac{\partial S_{\theta_2}^T}{\partial q}(p_1^{j}, q_0, \omega),
		\end{align}   
		for $P_1^{0}=p_0$ and $j=0,1,\dots,C$ such that  $\|p_1^{C+1} - p_1^{C} \|_2 \le \epsilon$
		for \(\epsilon: = 10^{-12}\). Then we obtain the system data at the next instant \((\widetilde{p}_1, \widetilde{q}_1)\) as follows:
		\begin{equation}
			\widetilde{p}_1 = p_1^{C},\quad \widetilde{q}_1 = q_0 + \frac{\partial S_{\theta_2}^T}{\partial p}(\widetilde{p}_1, q_0, \omega).
		\end{equation}
		
		\subsection{Training and prediction algorithm}
		The learning and prediction processes of the SGFNN model are detailed in Algorithm~\ref{alggfnnprocess}.

\begin{algorithm}[H]
	\caption{SGFNN model}
	\label{alggfnnprocess}  
	\begin{algorithmic}
		\State \textbf{Input:} $\mathcal{D}=\left\{\left(x_{0}^{(i)}, x_{1}^{(i)}\right)|i=1,\dots,M \right\} \,\,(M = N \times L)$. $x_0^{(i)}, x_1^{(i)} \in \mathbb{R}^{2d}$; Batch size $K$ and batch number $N_B$; Training set obtained from $\mathcal{D}$ using  `sub-sampling' strategy. 
		\State \textbf{Training:} Performed in batches
		
		\State Encoder part: 
		\begin{equation}
			B=\{z^{(i)} = E_\Delta(x_0^{(i)}, x_1^{(i)},\theta_1), \quad i = 1, \ldots, K\}
		\end{equation}
		
		\State 
		\begin{equation}
			\mathcal{L}_D(B) = \|\hat{f}_B - f_N\|_2 + \tau \cdot\left(\sum_{j=1}^6 \frac{1}{c_j}\left\|\widehat{\mu}_N^{(j)}(B)-\mu^{(j)}\left(P_{\mathcal{N}}\right)\right\|_2^2+\frac{\nu}{K} \sum_{1 \leq j<k \leq n_z}\left(\widehat{\rho}_{j k}(B)\right)^2\right)
		\end{equation}
		
		\State Decoder part: 
		\begin{equation}
			S_{\theta_2} = D_\Delta(p_1^{(i)}, q_0^{(i)}, z^{(i)},\theta_2), \quad (i = 1, \ldots, M)
		\end{equation}
		
		\State 
		\begin{equation}
			\mathcal{L}_{MSE}(B)=\frac{1}{K} \sum_{i=1}^K\left(\left\|p_1^{(i)}-p_0^{(i)}+\frac{\partial S_{\theta_2}^{T}}{\partial q_0^{(i)}}\right\|_2^2+\left\|q_1^{(i)}-q_0^{(i)}-\frac{\partial S_{\theta_2}^{T}}{\partial p_1^{(i)}}\right\|_2^2\right)
		\end{equation}
		
		\State Total loss:
		\begin{equation}
			\mathcal{L}_{SGFNN}=\frac{1}{N_B} \sum_{j=1}^{N_B} \mathcal{L}\left(B_j\right)=\frac{1}{N_B} \sum_{j=1}^{N_B}\left[\mathcal{L}_{MSE}\left(B_j\right)+\lambda \cdot \mathcal{L}_D\left(B_j\right)\right]
		\end{equation}
		
		\State \textbf{Output:} Minimizing $\mathcal{L}_{SGFNN}$ with {\bf Adam} to output $S_{\theta^*}$
		
		\State \textbf{Prediction:} Given $(p_0, q_0)$ and $\omega\sim\mathcal{N}(0,1)$, $(\widetilde{p}_1, \widetilde{q}_1)$ obtained by the symplectic scheme
		\begin{equation}
			\begin{aligned}
				\widetilde{p}_1 = p_0 - \frac{\partial S_{\theta_2^*}^T}{\partial q}(\widetilde{p}_1, q_0, \omega),\quad
				\widetilde{q}_1 = q_0 + \frac{\partial S_{\theta_2^*}^T}{\partial p}(\widetilde{p}_1, q_0, \omega).
			\end{aligned}
		\end{equation}
	\end{algorithmic}
\end{algorithm}
		
		\section{Numerical experiments}
		\label{sec4}
		
		In this section, we perform numerical experiments on our SGFNN model, and compare the numerical results with those of the sFML model (\cite{Xu}), demonstrating the superiority of our method. The experiments were conducted on a Windows 11 operating system with a virtual environment configured with Python 3.9.6 and TensorFlow 2.7.0.
		
		For each SHS, we uniformly use the midpoint method to generate 10,000 discrete trajectories as the observational data. The initial points are distributed across different regions, and the time step is set to \(\Delta = 0.01\). The lengths of the trajectories varies from system to system. The number of batches is $ n_B = 1000 $, and each batch contains $ K = 10000 $ data pairs. In the prediction stage, for each given initial value, we generate 10,000 predicted trajectories.
		
		The encoder and decoder have the same architecture, both with 3 hidden layers and 20 neurons per layer. The activation function used is ELU, and Adam is employed as the optimizer for training, with an initial learning rate of 0.001.
		
		\subsection{A linear stochastic oscillator}\label{4.1}
		Here we consider the linear stochastic oscillator with additive noise: 
		\begin{equation}\label{linear}
			\begin{aligned}
				& d p=-q d t+\sigma  d W, \quad p(0)=p_0, \\
				& d q=p d t, \quad q(0)=q_0,
			\end{aligned}
		\end{equation}
		which is a stochastic Hamiltonian system with Hamiltonians 
		\begin{equation}
			H_0 = \frac{1}{2}(p^2 + q^2),\quad
			H_1 = -\sigma q,
		\end{equation}
		where \(x = (p, q)^T \in \mathbb{R}^2\) and $\sigma\in\mathbb{R}$ is a number. It was proved that for the SHS \eqref{linear}, the second moment of the solution satisfies the linear growth property (\cite{higham})
		\begin{equation}\label{lineargr}
			E[p(t)^2+q(t)^2] = p_0^2+q_0^2 +\sigma^2 t.
		\end{equation}
		
		For data collection we sample the initial observational points \(\{(p_0^{(i)}, q_0^{(i)})|i=1,\dots,10000\}\) uniformly from the region \(U(0,3)\) representing the disc centered at the origin with radius 3, and the observational time span is till \(T = 10\). In the prediction stage, we take the initial value as \((p_0, q_0) = (0, 1)\), and the predicted time span is \(t\in [0,50]\).
		
		\begin{figure}[htbp]
			\centering
			\includegraphics[width=\linewidth,height=5cm]{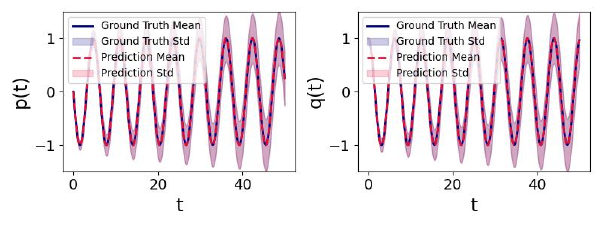}
			\label{f2.1}
			\includegraphics[width=\linewidth,height=5cm]{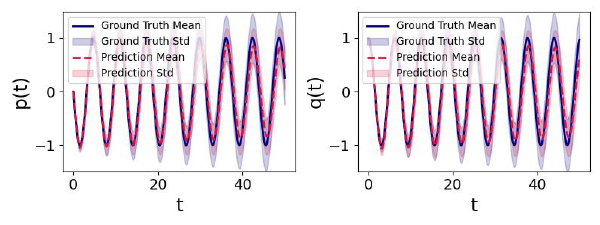}
			\label{f2.2}
			\caption{Linear Oscillator: comparison of the mean and standard deviation of 10,000 predicted trajectories and true trajectories, with \(\sigma = 0.1\). Upper: SGFNN model; Lower: sFML model.}
			\label{figlinearoscillatormeanstdsigma01}
		\end{figure}
		
		The upper and lower images in Figure \ref{figlinearoscillatormeanstdsigma01} present a comparison of the mean and standard deviation of the solutions generated by the SGFNN model and the sFML model with those of the true solution, respectively. The figures demonstrate that, for both the mean and the standard deviation, the solution predicted by the SGFNN model exhibits a higher degree of accuracy and closer alignment with the true solution compared to that predicted by the sFML model.
		
		\begin{figure}[htbp]
			\centering
			\includegraphics[width=7cm,height=5.5cm]{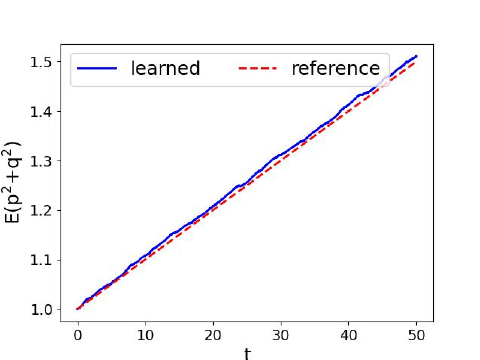}
			\label{f3.1}
			\includegraphics[width=7cm,height=5.5cm]{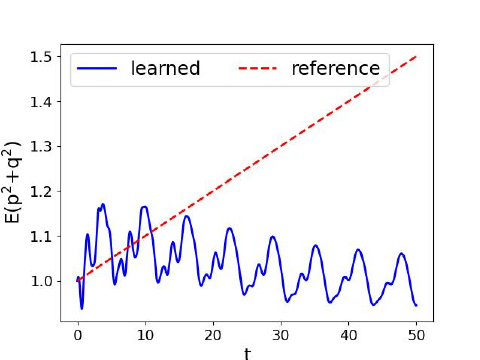}
			\label{f3.2}
			\caption{Linear Oscillator: growth of the second moment of the solution predicted by the SGFNN model (left) and the sFML model (right), with \(\sigma = 0.1\). }
			\label{figlinearoscillatorsecondmomentsigma01}
		\end{figure}
		Figure \ref{figlinearoscillatorsecondmomentsigma01} shows the evolution of the second moment $E(p(t)^2+q(t)^2)$ predicted by the SGFNN model (left) and the sFML model (right), respectively, on the time interval $[0,50]$, where the true evolution is plotted by red dashed lines as reference. It is clear that the SGFNN model can preserve the linear growth property of the second moment as described by \eqref{lineargr}, while the sFML model fails to maintain this property. 
		
		\begin{figure}[htbp]
			\centering
			\includegraphics[width=0.9\linewidth,height=5cm]{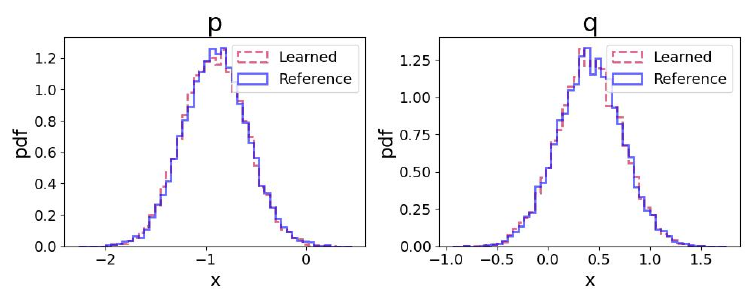}
			\label{f4.1}
			\includegraphics[width=0.9\linewidth,height=5cm]{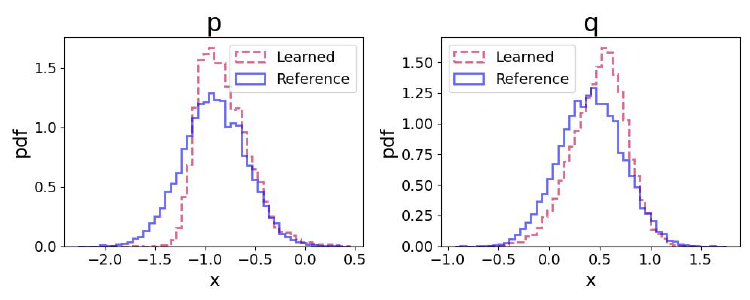}
			\label{f4.2}
			\caption{Linear Oscillator: histogram of the probability density function (pdf) of the solution predicted by the SGFNN model (upper) and the sFML model (lower) at \(T = 20\), with \(\sigma = 0.1\). }
			\label{figlinearoscillatorhistogramT1sigma01}
		\end{figure}
		Figure \ref{figlinearoscillatorhistogramT1sigma01} demonstrates the histogram of the probability density function (pdf) of the solution at $T=20$ predicted by the SGFNN model (upper) and the sFML model (lower), respectively, where that of the true solution is plotted by a blue line as reference. It can be seen that the SGFNN generates much more accurate pdf of the solution than the sFML model.
		
		\begin{figure}[htbp]
			\centering
			\includegraphics[width=0.5\linewidth,height=5cm]{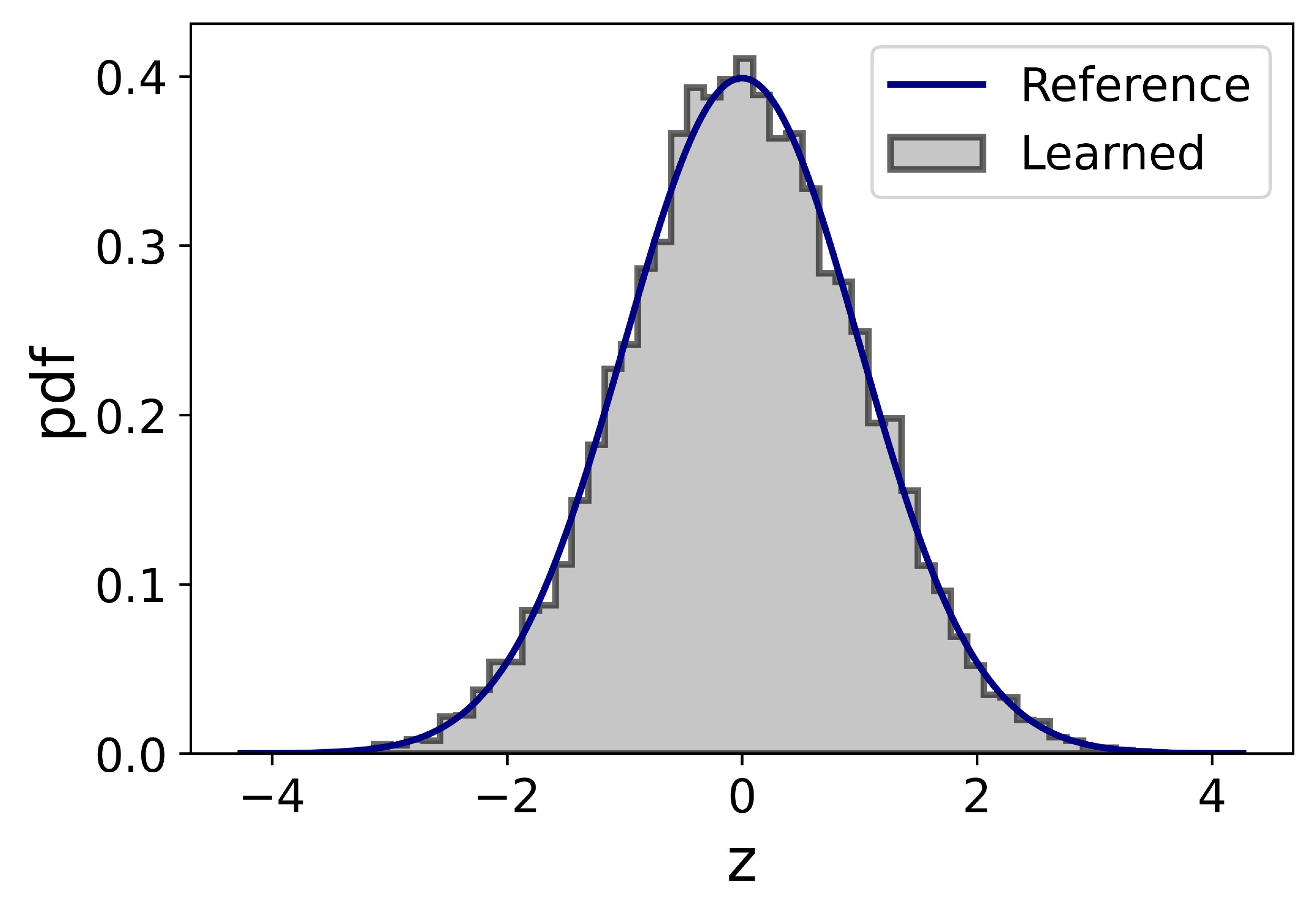} 
			\caption{Linear Oscillator: histogram of the pdf of the latent variable \(z\) generated by the encoder, with \(z = E_{\Delta}(x_0, x_1), x_0 = (0,1), \Delta=0.01, \sigma = 0.1\). Reference: Standard Normal Distribution.}
			\label{figlinearoscillatorlatentzsigma01}
		\end{figure}
		Figure \ref{figlinearoscillatorlatentzsigma01} shows the histogram of the pdf of the latent variable \(z\) generated by the trained encoder \(z = E_{\Delta}(x_0, x_1,\theta_1^*)\). Good agreement is observed between the trained pdf and the expected one . 
		
		\subsection{Kubo oscillator}
		In this example, we consider a stochastic Hamiltonian system with multiplicative noise, i.e. the Kubo oscillator:
		\begin{equation}
			\begin{aligned}
				&d p  =-a q d t-\sigma q \circ d W(t), & & p(0)=p_0, \\
				&d q  =a p d t+\sigma p \circ d W(t), & & q(0)=q_0,
			\end{aligned}
		\end{equation}
		where $(p,q)^T\in\mathbb{R}^2$, $a$ and $\sigma$ are constants, and the Hamiltonians are 
		\begin{equation}
			H_0 = \frac{1}{2} a (p^2 + q^2),\quad
			H_1 = \frac{1}{2} \sigma (p^2 + q^2).
		\end{equation}
		It is known that the quantity 
		$H(p, q) = p^2 + q^2$ for the Kubo oscillator is conservative (\cite{Mil_all}), namely
		\begin{equation}\label{circular}
			H(p(t),q(t))=p(t)^2+q(t)^2\equiv p_0^2+q_0^2,\quad \forall t\ge 0,
		\end{equation}
		that is, the phase flow of the Kubo oscillator is a circle centered at the origin with radius $\sqrt{p_0^2+q_0^2}$.
		
		Initial data \(\{(p_0^{(i)}, q_0^{(i)})|i=1,\dots,10000\}\) are uniformly sampled from the region \(U(0, 3)\), as in the example of section \ref{4.1}, and the observational time span is \(t\in[0,1]\). In the prediction stage, we take the initial value \(x_0=(p_0, q_0) = (1, 0)\), and the predicted time span is \(t\in[0,5]\).
		
		
		\begin{figure}[htbp]
			\centering 
			
			\includegraphics[width=\linewidth,height=5cm]{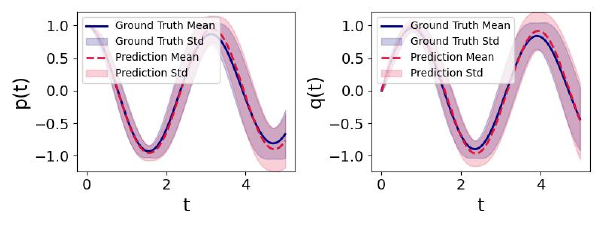}
			\label{f6.1} 
			
			\vspace{0.5cm}
			
			\includegraphics[width=\linewidth,height=5cm]{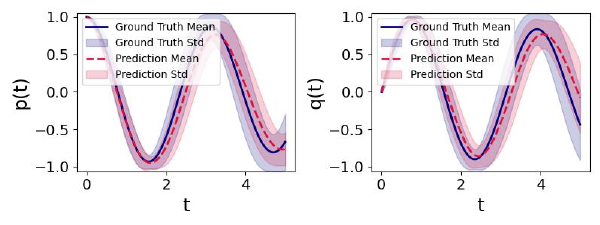}
			\label{f6.2} 
			
			\caption{Kubo Oscillator: comparison of the mean and standard deviation of 10,000 predicted trajectories and true trajectories, with \(a = 2\) and \(\sigma = 0.3\). Upper: SGFNN model; Lower: sFML model.}
			\label{figkubooscillatormeanstda2sigma03}
		\end{figure}

		Figure \ref{figkubooscillatormeanstda2sigma03} shows the comparison of the mean and standard deviation of the predicted solutions generated by the SGFNN model (upper) and the sFML model (lower) with those of the true solution, respectively. Obviously, the mean and standard deviation of the solution predicted by the SGFNN model deviate less from the true values compared to those predicted by the sFML model. 
		
		\begin{figure}[htbp]
			\centering
			\includegraphics[width=0.5\linewidth,height=6cm]{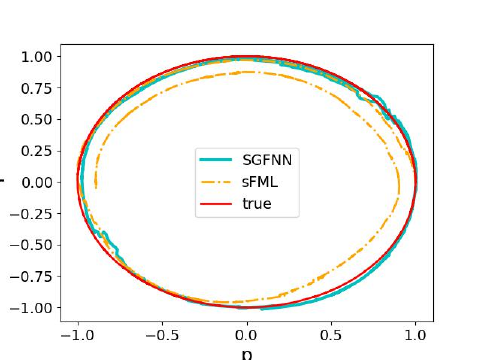} 
			\caption{Kubo Oscillator: predicted (p(t), q(t)) on $t\in[0, 5]$ by the SGFNN model and the sFML model, with \(a = 2\) and \(\sigma = 0.3\).}
			\label{figkubooscillatorpq} 
		\end{figure}
		
		Figure \ref{figkubooscillatorpq} illustrates a comparison of phase trajectories predicted by the SGFNN model and the sFML model with that of the true solution on $t\in[0,5]$, respectively. It is clear that the $p$-$q$ plots corresponding to SGFNN are closer to the real ones and preserve the circular property of the phase trajectory as described by \eqref{circular}, while the plots corresponding to the sFML model fails to preserve this property. 
		

		\begin{figure}[htbp]
			\centering 
			
			\includegraphics[width=0.9\linewidth,height=5cm]{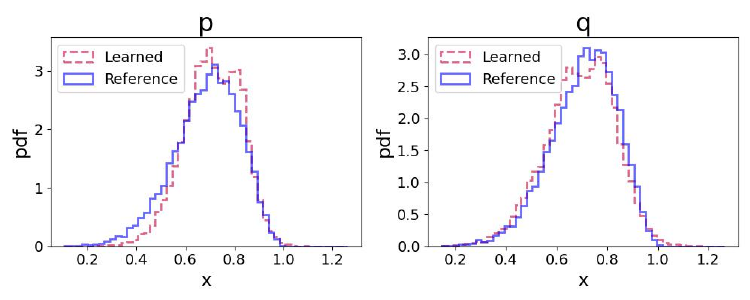}
			\label{f7.1} 
			
			\vspace{0.5cm}
			
			\includegraphics[width=0.9\linewidth,height=5cm]{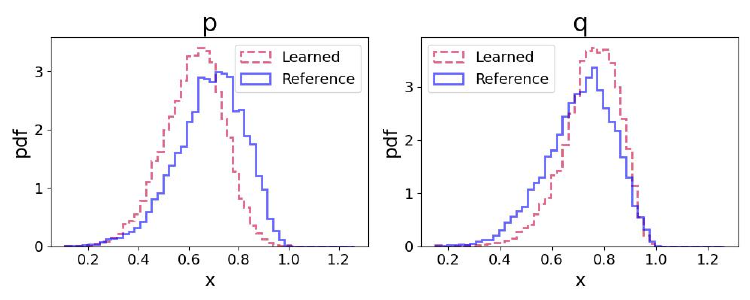}
			\label{f7.2} 
			
			\caption{Kubo Oscillator: histogram of the pdf of predicted and true solutions at \(T = 0.4\), with \(a = 2\) and \(\sigma = 0.3\). Upper: SGFNN model; Lower: sFML model.}
			\label{figkubooscillatorhistogramT04a2sigma03}
		\end{figure}
		
		Figure \ref{figkubooscillatorhistogramT04a2sigma03} presents the pdf histograms of the solutions predicted by the SGFNN model (upper) and the sFML model (lower) respectively at time \(T = 0.4\), where that of true solution is also plotted via blue lines as reference. It can be observed that the pdf predicted by the SGFNN model is closer to the true one than that by the sFML model.
		
		\begin{figure}
			\centering
			\includegraphics[width=0.5\textwidth,height=5cm]{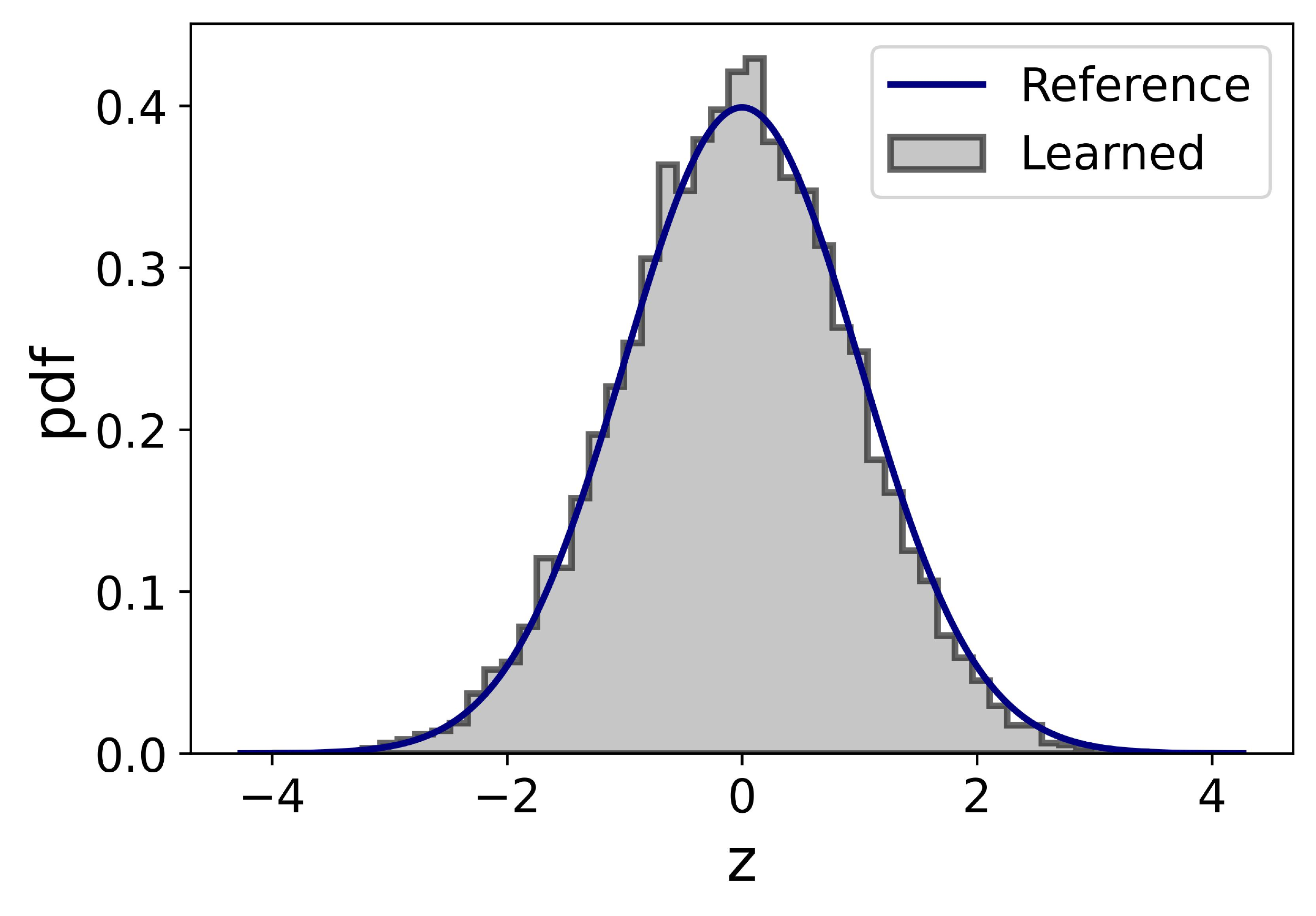} 
			\caption{Kubo Oscillator: histogram of the pdf of the latent variable \(z\) generated by the encoder, with \(a = 2\) and \(\sigma = 0.3\). Reference: Standard Normal Distribution.}
			\label{figkubooscillatorlatentza2sigma03} 
			\label{f8} 
		\end{figure}
		Figure \ref{figkubooscillatorlatentza2sigma03} illustrates that for the Kubo oscillator, the distribution of the latent variable \(z\) extracted by the encoder closely resembles a standard normal distribution.
		
		\subsection{A non-separable SHS}
		In this example, we consider a non-separable stochastic Hamiltonian system, where the Hamiltonians are
		\begin{equation}
			H_0 =\frac{1}{2}\left(p^2+1\right)\left(q^2+1\right), \quad 
			H_1 =0.1\left(p+q\right)^2
		\end{equation}
		with \(x = (p, q)^T \in \mathbb{R}^2\). 
		
		Initial points \(\{(p_0^{(i)}, q_0^{(i)}|i=1,\dots,10000)\}\) are uniformly sampled from the region \(U(0, 3)\), and the time span of the observational data is \(t\in [0,1]\).
		

		\begin{figure}[htbp]
			\centering 
			
			\includegraphics[width=\linewidth, height=5cm]{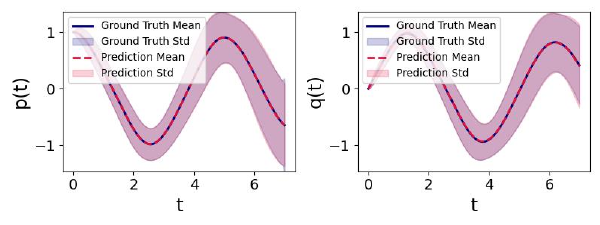}
			\label{f9.1} 
			
			\vspace{0.5cm}
			
			\includegraphics[width=\linewidth, height=5cm]{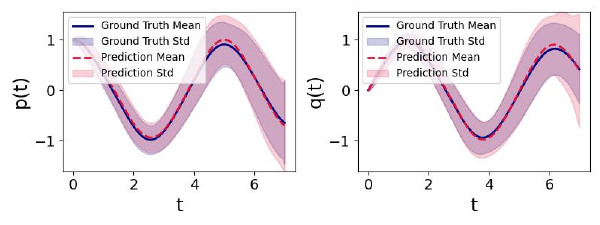}
			\label{f9.2} 
			
			\caption{Non-Separable SHS: comparison of the mean and standard deviation of 10,000 predicted trajectories and true trajectories. Top: SGFNN model; Bottom: sFML model.}
			\label{fignonseparableshsmeanstd}
		\end{figure}
		
		The complexity of the SHS increases in this scenario. However, as shown in Figure \ref{fignonseparableshsmeanstd}, the mean and standard deviation of the solution predicted by the SGFNN model still exhibit smaller deviations from the true ones compared to those predicted by the sFML model. 
		
		
		\begin{figure}[htbp]
			\centering 
			
			\includegraphics[width=0.9\linewidth, height=5cm]{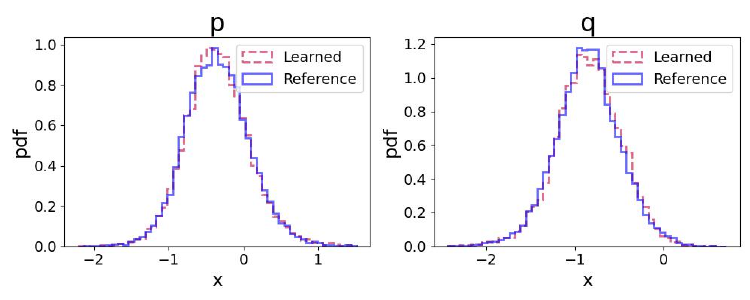}
			\label{f10.1} 
			
			\vspace{0.5cm}
			
			\includegraphics[width=0.9\linewidth, height=5cm]{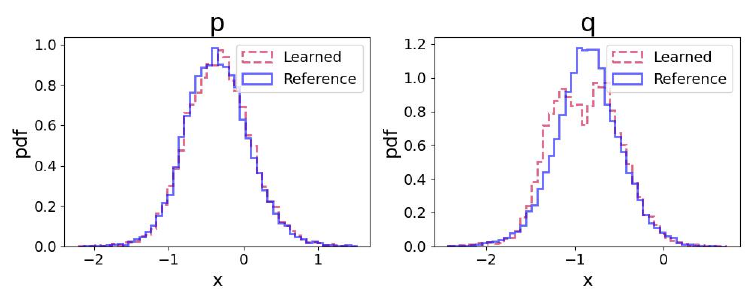}
			\label{f10.2} 
			
			\caption{Non-Separable SHS: histogram of the pdf of predicted and true solutions at \(T = 3.5\). Top: SGFNN model; Bottom: sFML model.}
			\label{fignonseparableshshistogramT35}
		\end{figure}
		
		Figure \ref{fignonseparableshshistogramT35} displays the distribution of the predicted solutions at time \(T = 3.5\) by the SGFNN (upper) and sFML models (lower), respectively. Similar to previous findings, the performance of the SGFNN model surpasses that of the sFML model, demonstrating a closer match to the true distribution.
		
		\begin{figure}
			\centering
			\includegraphics[width=0.5\textwidth,height=5cm]{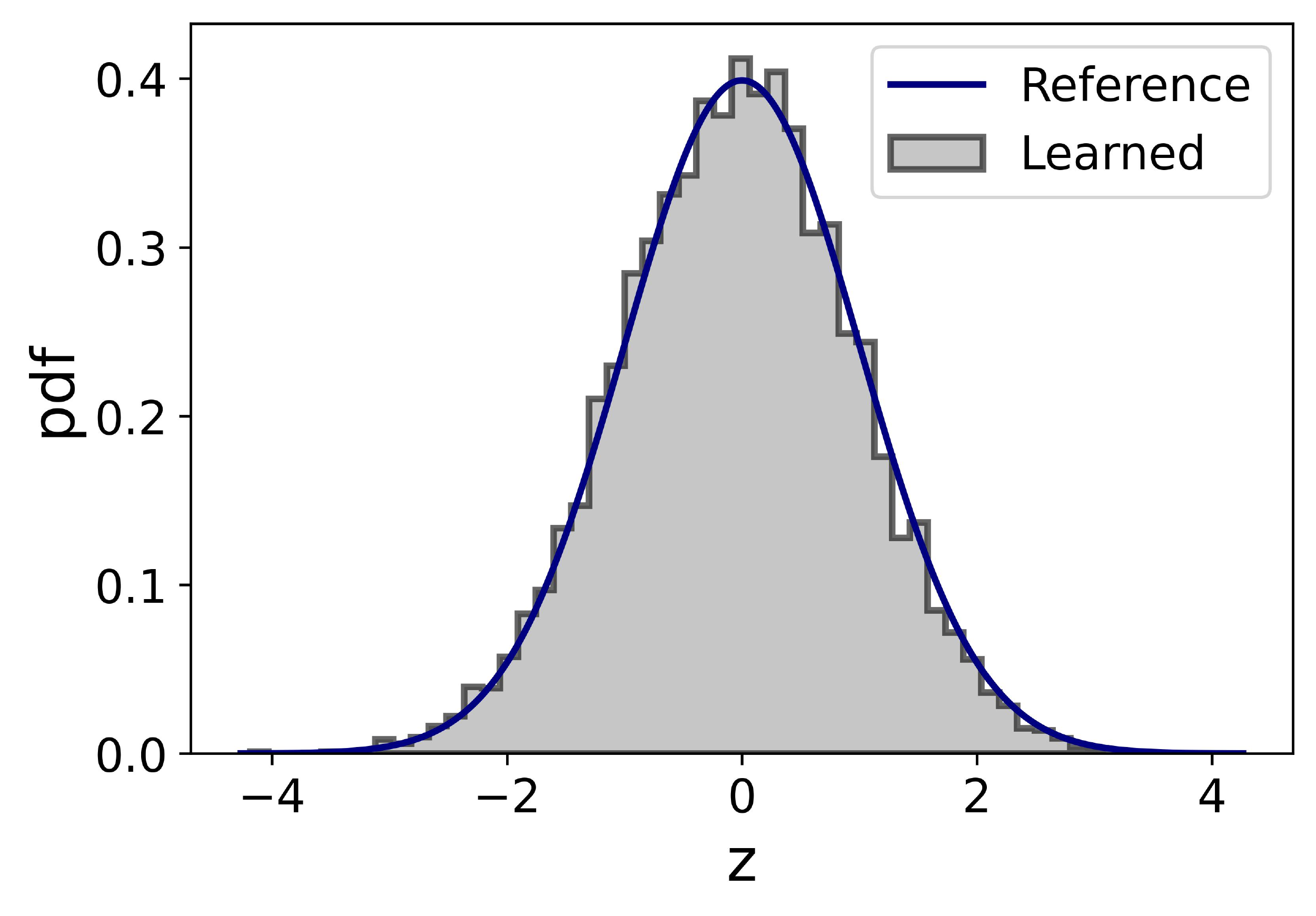} 
			\caption{Non-Separable SHS: histogram of the pdf of the latent variable \(z\) generated by the encoder. Reference: Standard Normal Distribution.}
			\label{fignonseparablshslatentz} 
		\end{figure}
		Figure \ref{fignonseparablshslatentz} illustrates the distribution of the latent variable \(z\) extracted by the encoder, which well resembles the standard normal distribution as desired, showing effectiveness of the SGFNN model in simulating the latent stochasticity even for complex systems.
		
		\subsection{Synchrotron oscillations: a SHS with two multiplicative noises}
		In this section, we consider the model of synchrotron oscillations which is described by the following SHS with two multiplicative noises (\cite{Mil_all}):
		
		\begin{equation}
			\begin{aligned}
				&dp = -\omega^2\sin(q) \, dt - \sigma_1\cos(q) \circ dW_1 - \sigma_2\sin(q) \circ dW_2, \quad p(0) = p_0, \\
				&dq = p \, dt, \quad q(0) = q_0,
			\end{aligned}
		\end{equation}
		where $(p, q) \in \mathbb{R}^2$, $\omega$, $\sigma_1$, $\sigma_2$ are constants, and the Hamiltonians are 
		\begin{equation}
			H_0 = -\omega^2\cos(q) + \frac{p^2}{2}, \quad H_1 = \sigma_1\sin(q), \quad H_2 = -\sigma_2\cos(q).
		\end{equation}
		
		
		\begin{figure}[htbp]
			\centering 
			
			\includegraphics[width=\linewidth, height=5cm]{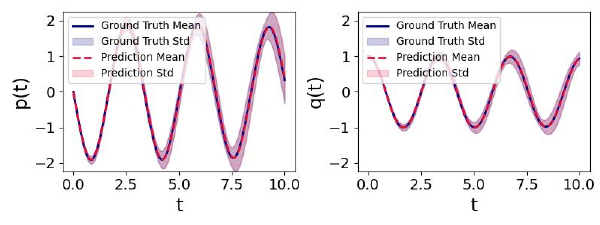}
			\label{f11.1}
			
			\vspace{0.5cm}
			
			\includegraphics[width=\linewidth, height=5cm]{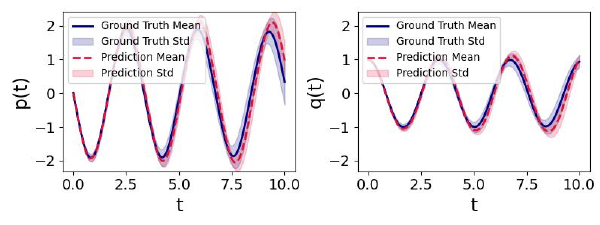}
			\label{f11.2}
			
			\caption{Synchrotron Oscillations: comparison of the mean and standard deviation of 10,000 predicted trajectories and true trajectories. Upper: SGFNN model; Lower: sFML model.}
			\label{figsynchrotron}
		\end{figure}
		
		Figure \ref{figsynchrotron} shows the comparison of the mean and variance of 10,000 trajectories predicted by the SGFNN model (upper) and the sFML model (lower) with those of the true trajectories, respectively. It is evident that the moments predicted by the SGFNN model are closer to that of the true solutions, and the prediction bias of the sFML model increases with time. 
		

		\begin{figure}[htbp]
			\centering 
			
			\includegraphics[width=0.9\linewidth, height=5cm]{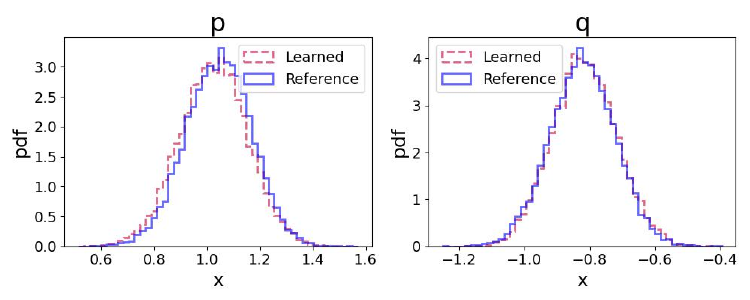}
			\label{f12.1}
			
			\vspace{0.5cm}
			
			\includegraphics[width=0.9\linewidth, height=5cm]{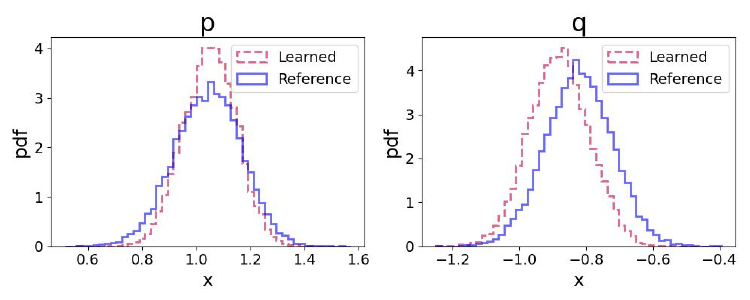}
			\label{f12.2}
			
			\caption{Synchrotron Oscillations: histogram of the pdf of predicted and true solutions at \(T = 2\). Upper: SGFNN model; Lower: sFML model.}
			\label{figsynchrotronhistogramT35}
		\end{figure}
		
		We also compare the distribution histograms of the predicted solutions generated by both models at a specific time point $T=2$, as shown in the Figure \ref{figsynchrotronhistogramT35}. Again, the results demonstrate that the SGFNN model outperforms the sFML model.
		
		\begin{figure}[htbp]
			\centering
			\includegraphics[width=0.45\linewidth,height=5cm]{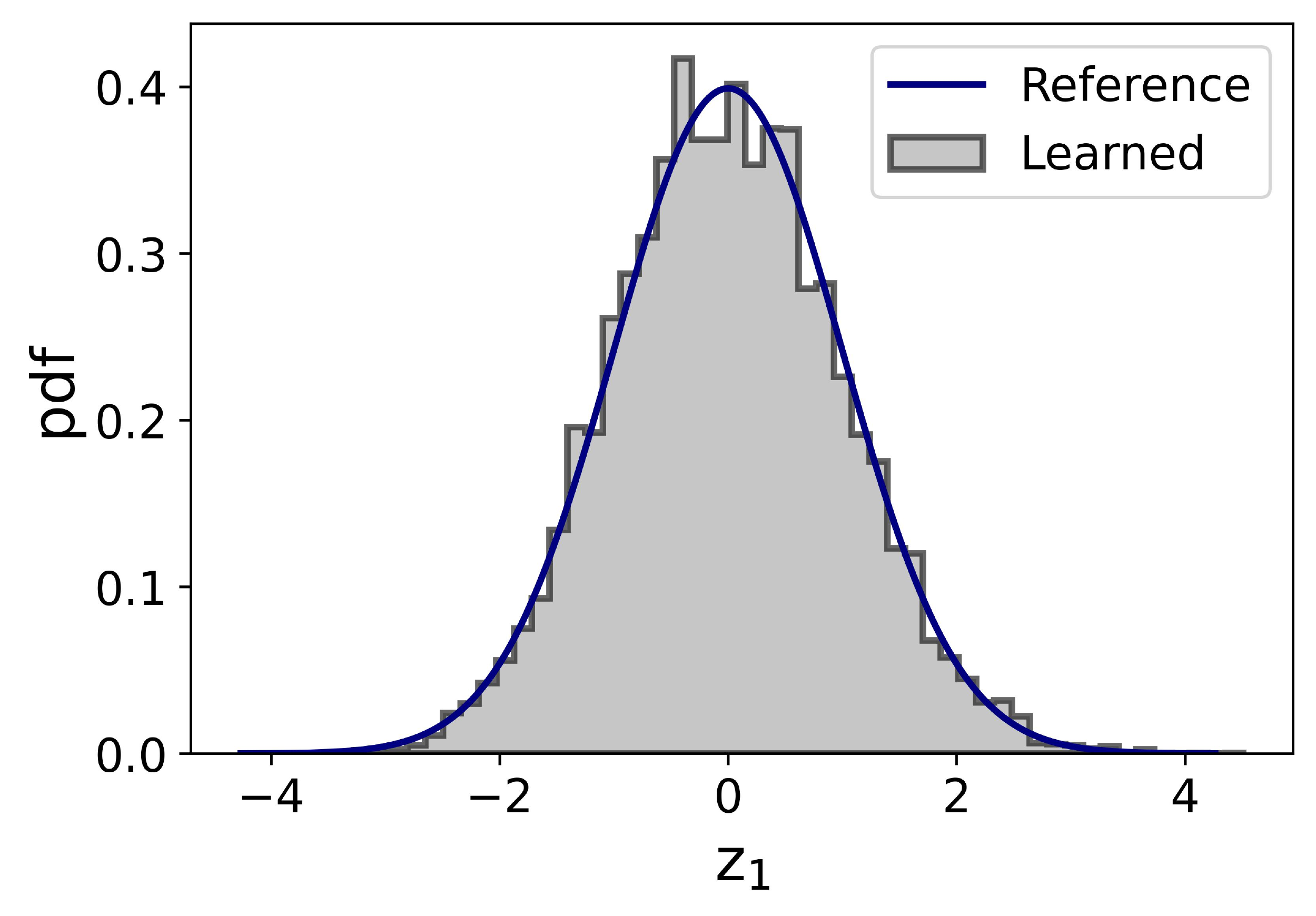}
			\label{f13.1}
			\includegraphics[width=0.45\linewidth,height=5cm]{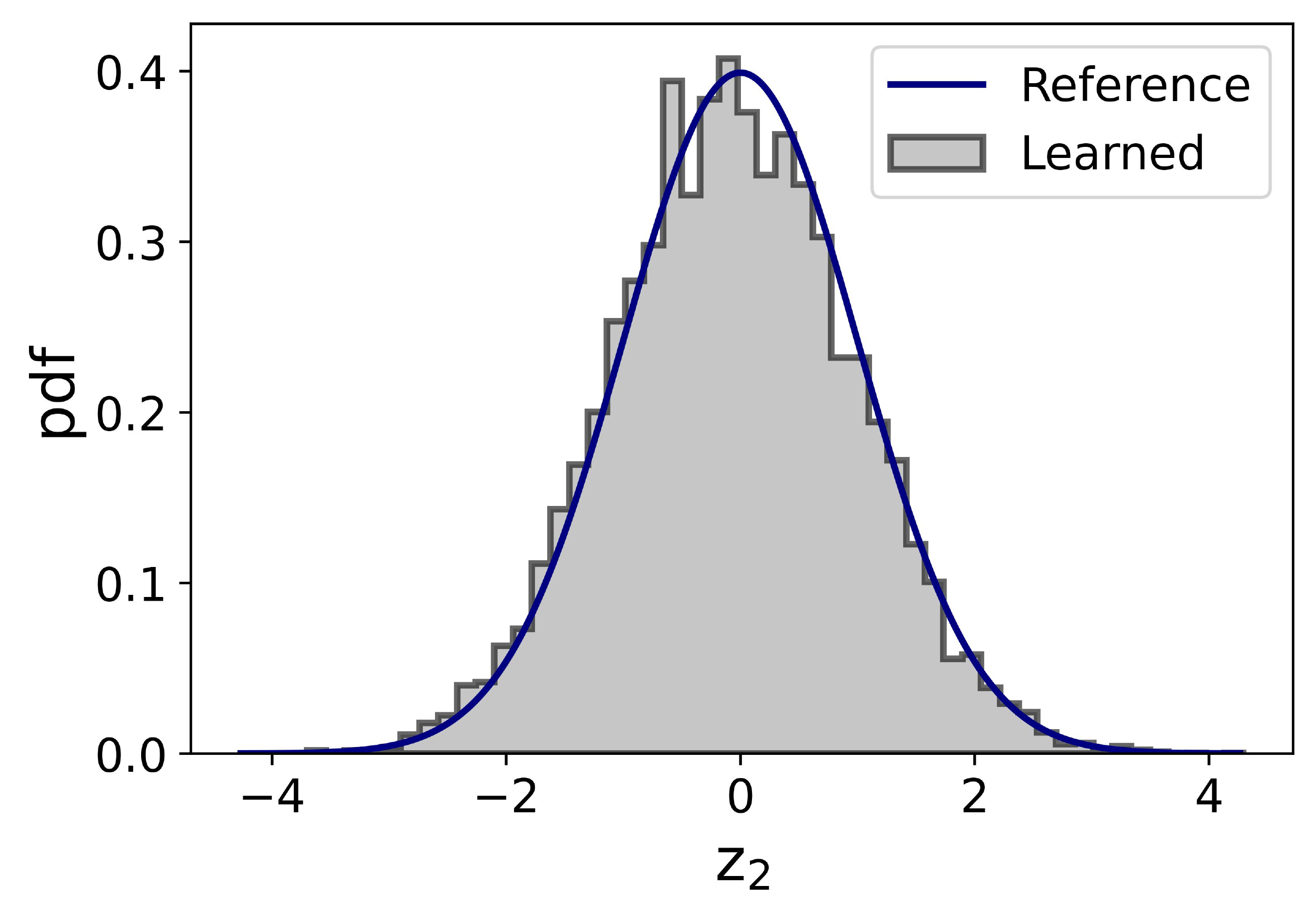}
			\label{f13.2}
			\caption{Synchrotron Oscillations: histogram of the pdf of the latent variables \(z_1\) and $z_2$ generated by the encoder. Reference: Standard Normal Distribution.}
			\label{figsynchrotronplotz}
		\end{figure}
		
		For the synchrotron oscillation system, the latent variable $z=(z_1,z_2)$ extracted by encoder should follow a two-dimensional unit Gaussian distribution. Figure~\ref{figsynchrotronplotz} shows the marginal density functions of $z_1$ and $z_2$ generated by the encoder respectively. It can be observed that the distributions of both $z_1$ and $z_2$ closely resemble the standard normal distribution.
		
		Beyond the figures presented above, we also compare 
		the two models via calculating the following quantities as in \cite{Liu1}, across the four examples mentioned above:
		
		\begin{itemize}
			\item Difference in mean between the predicted $\widetilde{x}_T$ and true values $x_T$ at the end of the prediction time: $e_{T}^{m}=\|E(x_{T})-E(\widetilde{x}_{T})\|$;
			\item Difference in standard deviation between the predicted and true values at the end of the prediction time: $e_{T}^{std}=\|STD(x_{T})-STD(\widetilde{x}_{T})\|$.
		\end{itemize}
		The results are summarized in Table~\ref{table}. Better results are highlighted in bold. It is clear that all the compared errors of the SGFNN model are smaller than that of the sFML model. This powerfully demonstrates the superiority of the structure-preserving learning which incorporates the physical characteristics into the neural network architecture, compared to general learning methods.
		
		\begin{table}[H]
			\caption{\textbf{Comparison of mean and STD estimations by SGFNN and sFML models}}
			\label{table}
			\centering
			\begin{tabular}{cccccc}
				\toprule
				Example &T & $e^m_T$ of sFML & $e^m_T$ of SGFNN &$e^{std}_T$ of sFML & $e^{std}_T$ of SGFNN \\
				\midrule
				Linear\_0.1&50&3.5349E-01&\textbf{2.8755E-02}&1.3608E-01&\textbf{3.5269E-03} \\
				Kubo&5&3.6795E-01&\textbf{9.8468E-02}&1.5113E-01&\textbf{1.2450E-01} \\
				non-separable SHS&7&6.3663E-02&\textbf{1.4602E-02}&4.4735E-01&\textbf{7.8563E-02}\\
				Synchrotron Oscillations&10&1.0791E-00&\textbf{7.4355E-02}&3.7287E-01&\textbf{1.1828E-02}\\
				\bottomrule
			\end{tabular}
		\end{table}
		
		\section{Conclusion}\label{sec5}
		We put forward a neural network model for learning the unknown stochastic Hamiltonian system from observed data, employing an autoencoder network model to acquire the system's stochastic generating function, which can generate the symplectic flow map almost surely for structure-preserving prediction of the underlying stochastic Hamiltonian system. The encoder and decoder parts play the role of capturing the random variable hidden inside data pairs, and the stochastic generating function of the system, respectively, via imposed fitting to desired distribution and observed data by the loss function. 
		Numerical experiments compare our method with the benchmark stochastic flow map learning (sFML) model on four specific stochastic Hamiltonian systems. Our model exhibits better learning and prediction accuracy while possessing the capacity of maintaining the system's symplectic structure. The results indicate that respecting and preserving the intrinsic structure of the underlying dynamical systems in the neural network architecture can offer significant improvement in the learning behavior of neural networks.
		
		\section*{Acknowledgments}
		The first, second and fourth authors are supported by the NNSFC No. 11971458.


\end{document}